\documentclass[12pt]{article}

\usepackage{amssymb}
\usepackage{amsmath}
\usepackage{latexsym}
\usepackage{amsfonts}
\usepackage[latin1]{inputenc}
\usepackage{graphicx}
\usepackage{subfigure}
\usepackage{epsfig}
\usepackage{color}
\usepackage{dsfont}
\usepackage{algorithm,algorithmic}

\usepackage[title]{appendix}

\newcommand{\bfr}{{\bf{r}}}
\newcommand{\bfv}{{\bf{v}}}
\newcommand{\bfu}{{\bf{u}}}
\newcommand{\pos}{\mathbf{x}}
\newcommand{\node}{\mathbf{x}}
\newcommand{\normal}{\mathbf{n}}
\newcommand{\nn}{\mathbf{n}}

\newcommand{\R}{\mathbb{R}}
\newcommand{\Sph}{\mathbb{S}^2}

\newcommand{\ds}{\displaystyle}
\newcommand{\nablav}{{\nabla}}
\newcommand{\nablat}{\nablav_{\Sigma}}

\begin{document}

%\begin{frontmatter}
\title{A local radial basis function method for  the Laplace-Beltrami operator}
\author{Diego \'Alvarez, Pedro Gonzalez-Rodriguez, and Manuel Kindelan}
%\author{Pedro Gonzalez-Rodriguez\corauthref{label1}}, and 
%\author{Manuel Kindelan}, 
%% \author{Miguel Moscoso\corauthref{label1}}
%\ead{jdiego@math.uc3m.es}
%\ead{pgonzale@ing.uc3m.es}
%\ead{kinde@ing.uc3m.es}
%% \ead{moscoso@math.uc3m.es}\corauth[label1]{Corresponding author.
%% Address: Universidad Carlos III de Madrid, Avenida de la
%% Universidad 30, 28911 Legan\'es, Spain. Fax: +34 91 624 91 29}
%\author{}
%\address{Universidad Carlos III de Madrid, Departamento de Matem\'aticas, 
%Avenida de la Universidad 30, 28911 Legan\'es, Spain}

\maketitle
\begin{abstract}
We introduce a new local meshfree method for the approximation of the Laplace-Beltrami operator on a smooth surface of co-dimension one embedded in $\R^3$. A key element of this method is that it does not need an explicit expression of the surface, which can be simply defined by a set of scattered nodes. It does not require expressions for the surface normal vectors and for the curvature of the surface neither, which are approximated  using formulas derived in the paper. An additional advantage is that it is  a local method and, hence, the matrix that approximates the Laplace-Beltrami operator is sparse, which translates into good scalability properties. 
The convergence, accuracy and other computational characteristics of the method are studied numerically. The performance 
is shown by solving two reaction-diffusion partial differential equations on surfaces; the Turing model for pattern formation, and the Schaeffer's model for electrical cardiac tissue behavior.

\end{abstract}
%s
%\begin{keyword}
%Node selection; finite differences
%
%%\PACS 02.30.Zz, 42.30.Wb, 41.20.Cv_
%\end{keyword}
%\end{frontmatter}

\section{Introduction}

Researchers from a diverse range of areas such as medicine, geoscience or computational graphics  \cite{Alvarez12,Barreira11,Carr01,Flyer11,Gu04} 
often urge to find solutions to partial differential equations (PDEs) on surfaces. This is a challenging problem due to the complexity of the differential operators 
restricted to the geometry of the given surfaces. Among all the methods that seek the solution of this class of PDEs, the most popular are finite element methods, where the equations are solved using surface triangulation \cite{Dziuk07,Holst01,Turk91}. They are in general very efficient but they also have several problems. Specifically, the discretization may not be trivial, and some difficulties may appear when computing  geometric primitives such as surface normals and curvatures. Another common approach is to embed the surface PDE within a differential equation posed on the whole $\R^3$, and then restrict the solution to the surface of interest \cite{Bertalmio01,Greer06,Macdonald09}. The main advantage of this method comes from working with cartesian grids. However, the discretization is done in a higher dimensional space and, thus, the computational cost increases.
Another disadvantage of these methods is that  it is not well understood how the accuracy deteriorates in this procedure.

An alternative approach to estimating the solutions to PDEs on surfaces are  global RBF-based methods \cite{Flyer09,Fuselier13,Lehto17}.  They can naturally handle irregular geometries and scattered nodes layouts as no triangularization is needed. Their main advantages are the ease of implementation and the potential spectral accuracy with respect to the number of nodes $n$ used to discretize the surface. An additional advantage is that they operate using cartesian coordinates instead of intrinsic coordinates over the surface. However, there are at least two important drawbacks. One is that the computational cost scales as $O(n^3)$, and thus, it becomes rapidly excessive, making global RBF methods impractical for large problems. The other drawback is that there exists an inverse correlation between the accuracy and the stability of these methods, being the shape parameter that which determines the trade-off between the two. Indeed, they need a somewhat ad-hoc choice of the value of a shape parameter that affects the stability and the ill-conditioning of the resulting linear systems.

To overcome these drawbacks local versions of the global RBF method have also been proposed to solve PDEs defined on a surface \cite{Piret2012,Shankar2013,Shankar15}. These  methods use only a relatively small subset of all the available points to approximate the PDE operator locally, so their cost is largely reduced and the scaling properties improved \cite{fornberg2011}. 
As a side effect, the spectral accuracy is sacrificed but much better conditioned linear systems are obtain. In the fewest words possible, 
local RBF methods inherit many  key strengths of global RBF method, with a reduced cost but loosing accuracy.

An interesting consequence of local RBF methods is that they generate finite differences (FD) schemes with certain weights $w_i$.
Hence, they are known as RBF-FD methods. They benefit from some of the key properties of traditional FD approximations, but they are more flexible. 
While FD are enforced to be exact for polynomials evaluated at the node $x_k$, RBF-FD are enforced to be exact for RBF interpolants. Thus,
weights for FD-like stencils with  scattered nodes can be obtained easily.

%As for FD, these formulas are generated using only a local set of nodes (stencils) so that the resulting differential matrices are sparse like in the standard FD method. This sparsity makes RBF-FD an attractive alternative to global RBF methods as they perform better in terms of accuracy per computational cost \cite{fornberg2011}. Additionally, the local character of RBF-FD methods makes them more flexible in terms of local grid refinement strategies than global RBF methods. The main drawback of local methods is that they only have algebraic convergence. 

%Fuselier and Wright \cite{Fuselier13} used the RBF global method to solve reaction-diffusion equations on smooth, closed embedded surfaces. They used RBF interpolants to approximate the Laplace-Beltrami operator at a set of scattered nodes on a given surface. Their method combines the advantages of intrinsic methods with those of embedded methods. More recently Shankar {\em et al.} \cite{Shankar15} proposed a similar method based on the local RBF-FD approach. Their method requires only scattered nodes representing the surface and normal vectors at those scattered nodes. A more accurate version of the method is based on Hermite RBF interpolation \cite{Shankar15}.

In  \cite{Alvarez18}, we derived a closed-form formula for the (global)  RBF-based approximation of the Laplace Beltrami Operator (LBO). 
It can not only be applied to surfaces whose explicit formula is known, but also to surfaces defined by a set of scattered nodes. The formula only
requieres knowledge on the positions of those points, and its accuracy  mainly depends on the errors made in the RBF surface reconstruction.
If an explicit representation of the surface is given and, therefore, the normal vectors to the surface and its curvature are computed analytically, the convergence  
is faster than algebraic. 

In this paper, we generalize the approach in  \cite{Alvarez18} to local RBF approximations. We carry out a thorough study of the proposed
method and analyze its convergence and stability properties when it is applied to reaction diffusion equations over surfaces. We consider
the Turing model for pattern formation in nature, and the Schaeffer's model for electrical cardiac tissue behavior.

%This involves the approximation of the superficial gradient and the Laplace-Beltrami operator (LBO). This method uses Cartesian coordinates, thereby avoiding the singularities typically associated with intrinsic coordinate systems. 
%The method has less computational cost than the global RBF method  while still retaining the ability to use scattered points on the surface to approximate derivatives, combining the benefits of the intrinsic and narrow-band approach. Another important strength of our method is that we do not need an implicitly defined surface, which is a big advantage over previous methods.

The paper is organized as follows. In Section 2, 
%we give the formulation of the problem, 
%starting with a short explanation of the RBF-FD method. Then 
we derive an analytical expression of the LBO applied to a generic RBF,
%which will be used to calculate the weights of the RBF-FD approximation to this operator. Next 
and we characterize the surface defined by a cloud of points in order to approximate the normal vectors and the curvature. In Section 3, we 
analyze the convergence and stability properties of the proposed method.
%present some numerical experiments that show the performance of the approximation to the LBO and to the characterization of the surface. We have also carried out a stability analysis of this approximation calculating its eigenvalues and comparing them with the analytical ones. 
In Section 4, we  apply our approximation to the solution of two different PDEs over surfaces. 
%The first one is the solution of the Turing model for pattern formation in nature, and the second one is the Schaeffer's model for the cardiac tissue behavior. 
Section 5 contains our conclusions.
%We finish the paper with the conclusions and acknowledgments. 

\section{Formulation of the problem}
\label{sec:formulation}

Let $\Sigma$ be a two dimensional surface embedded in $\R^3$, and   ${\bf{X}}=\{\pos_i\}_{i=1}^n$  a set of $n$ scattered points  distributed on $\Sigma$. The objective of this section is to %obtain 
derive a local RBF-FD approximation 
%type   local RBF-FD
for the  Laplace-Beltrami operator (LBO) applied to a scalar function $f:\Sigma\to\R$  using the data  $(\pos_i, f(\pos_i))$.  We start reviewing the basic concepts of  RBF methods necessary to calculate the numerical approximation of a differential operator; see \cite{Buhmann03,Sarra09,Fonberg15} for more details. %After this we proceed to the calculation 
Then, we describe how to compute  the weights for local RBF-FD approximation of  the LBO. As in \cite{Alvarez18},  it is also necessary to compute some characteristics of the surface, such as the normal vectors and its curvature, which are also %be calculated 
computed using a local RBF methodology. %The final objective will be to build a 
Finally, we describe how to build a differentiation matrix that %allows to incorporate 
approximates the action of the LBO in reaction-diffusion type equations on surfaces.

\subsection{RBF fundamentals}
\label{subsec:formulation}
Radial functions are real valued functions $\Phi:{\mathbb {R}}^3\to \R$, whose value depends only on the magnitude of its argument, i.e., $\Phi(\pos)=\phi(\Vert\pos\Vert) \,,\,\pos \in{\mathbb {R}}^3\ $, where $\|\cdot \|$ is, usually, the standard euclidean norm. These functions $\phi:[0,\infty)\to\R$ are called  Radial Basis Functions (RBF). 
RBFs can be used to construct  interpolants for continuous functions $f:\Omega \to \R$,  $\Omega \subseteq \R^3$, sampled at a set of points $\{\pos_i\}_{i=1}^n\in\Omega$.  Given a set of 
data $\{\pos_i,f_i\}_{i=1}^n$,%for $i=1,2,\cdots,n$,  
the interpolant 
\begin{equation}
\label{eq:interpt}
S_f(\pos)=\sum_{i=1}^n \lambda_i\phi(\Vert\pos-\pos_i\Vert),\quad \pos\in\Omega,
\end{equation} 
is constructed as a linear combination of translates of $\phi$ over the points $\pos_i$, 
where the interpolation coefficients $\lambda_i$ are calculated %imposing 
by imposing that the interpolant coincides with the function at each point $\pos_i$, so
\begin{equation}
\label{eq:cond1}
S_f(\pos_i)=f(\pos_i),\,\hspace{1cm} i=1,2,...,n.
\end{equation}
These requirements lead to the following linear system
\begin{equation}\label{eq:interpolacion}
\underbrace{
	\begin{bmatrix}
	\phi(r_1(\pos_1)) & \phi(r_1(\pos_2)) &  \hdots  & \phi(r_1(\pos_n))  \\ 
	\phi(r_{2}(\pos_1)) & \phi(r_2(\pos_2)) &  \hdots  & \phi(r_{2}(\pos_n)) \\ 
	\vdots  &  \vdots  &\ddots  & \vdots     \\ 
	\phi(r_{n}(\pos_1)) & \phi(r_{n}(\pos_2)) &  \hdots  & \phi(r_n(\pos_n))  \\ 
	\end{bmatrix}}_{{{A}}}
\underbrace{
	\begin{bmatrix}
	\lambda_1 \\ 
	\lambda_2 \\ 
	\vdots  \\ 
	\lambda_n \\ 
	\end{bmatrix}
}_{\boldsymbol{\lambda}}
=\underbrace{\begin{bmatrix}
	f(\pos_1) \\
	f(\pos_2)  \\ 
	\vdots  \\ 
	f(\pos_n) \\ 
	\end{bmatrix}}_{\bf{f}},
\end{equation}
where $r_i(\pos)=\lVert\bfr_i(\pos)\rVert=\lVert \pos- \pos_i \rVert$. 
Notice that the RBF interpolation matrix $A$ is dense, symmetric and, under certain conditions, invertible \cite{Buhmann03}.
%Solving this linear system we obtain the interpoaltion coefficients 

It is sometimes useful to add 
a constant %constants 
to the RBF interpolant \eqref{eq:interpt}, so %in the form 
\begin{equation}
\label{eq:interpt2}
S_f(\pos)=\sum_{i=1}^n \lambda_i\phi(r_i(\pos))+\gamma\, .
%\sum_{j=1}^m \gamma_j p_i(\pos),\quad \pos\in\R^{m},
\end{equation} 
%where $\{p_j\}_{j=1}^{m}:\R^m\to\R$ is a basis of linear space of polynomial of degree less than $m$, 
Adding a
constant to the RBF interpolant  %allows 
ensures the exact interpolation of constant functions. % to obtain 
It also results in a less oscillatory interpolant, %improving 
and improves the accuracy of the interpolation \cite{Fonberg15}.
To compute the interpolation coefficients $\lambda_i$ and the constant $\gamma$ 
we impose the conditions (\ref{eq:cond1}) 
with the additional constraint
\begin{equation}
\label{eq:ligadura2}
%\sum_{i=1}^{n}\lambda_i p_j(\pos_i)=0, \quad j=1,2,\cdots,m.
\sum_{i=1}^{n}\lambda_i=0.%, \quad j=1,2,\cdots,m.
\end{equation}
Then, equation \eqref{eq:interpolacion} takes the form
\begin{equation}\label{eq:interpolacion_c}
\underbrace{
	\begin{bmatrix}
	A& {\bf{c}}  \\ 
	{\bf{c}}^T& 0 
	\end{bmatrix}}_{{{A_c}}}
\underbrace{
	\begin{bmatrix}
	{\boldsymbol{\lambda}}\\ 
	{{\gamma}}
	\end{bmatrix}
}_{\boldsymbol{\lambda}_c}
=\underbrace{\begin{bmatrix}
	{\bf{f}}\\
	{{0}} 
	\end{bmatrix}}_{{\bf{f}}_c},
\end{equation}
where $A$, ${\boldsymbol{\lambda}}$  and ${\bf{f}}$ are defined in  \eqref{eq:interpolacion}, and ${\bf{c}}$ is  a vector with $n$ components all equal to 1.
Thus, the interpolation problem can be solved calculating the coefficients as
\begin{eqnarray}\label{eq:lambdac}
{\boldsymbol{\lambda}}_c= A_c^{-1}{\bf{f}}_c.
\end{eqnarray}

RBFs can also be used to approximate differential operators %applied to a function 
in a similar way %that it 
to what is done in the FD method.
FD formulas approximate differential operators applied to a function $f$, at a point $\pos_0$, by a weighted sum 
\begin{eqnarray}\label{eq:weight}
\mathcal{L}_f{(\pos_0)}\approx \sum_{i=1}^n w_i\, f{(\pos_i)}
\end{eqnarray}
of the values of that function at a set $\{\pos_i\}_{i=1}^n$ of $n$ neighboring points of $\pos$.
In the standard FD formulation %, these weights,  
the weights $w_i$ are computed by enforcing (\ref{eq:weight}) to be exact for polynomials of a certain degree evaluated at the set of points $\{\pos_i\}_{i=1}^n$, while for RBF, the formula is enforced to be exact for RBF interpolants 
at those same points \cite{Fonberg15}. The resulting formulas are called Radial Basis Function-Finite Difference (RBF-FD) formulas.

Thus, given a differential operator $\mathcal{L}$ and a set of points, %to calculate 
the RBF-FD  weights are computed %we impose 
by imposing (\ref{eq:weight}) to be exact for the interpolant (\ref{eq:interpt2}), so
\begin{eqnarray}\label{eq:condition}
\sum_{i=1}^{n} w_i f(\pos_i)=\mathcal{L}_{S_f}(\pos_0)=\sum_{i=1}^n \lambda_i \mathcal{L}_{\phi(r_i)}(\pos_0)+ \gamma\mathcal{L}_1. %\sum_{j=1}^m \gamma_{j}\mathcal{L}_{p_j}(\pos).%\quad i=1\cdots n.
\end{eqnarray}
Equation (\ref{eq:condition}) can be written as
% 
%If we define the vector $\boldsymbol{ {L}}_{c} (\pos)=[ \mathcal{L}_{\phi(r_1)},\cdots, \mathcal{L}_{\phi(r_n)}(\pos), \mathcal{L}_1]^t$  and use \eqref{eq:lambdac}, 
%%and the symmetry of $A$,
% this equation can be written  as
\begin{eqnarray}\label{eq:Vweight}
%{\bf{f}}^t\,{\bf{w}}={\boldsymbol{\lambda}}^t\,\boldsymbol{ {L}}_{p} (\pos)= {\bf{f}}^t \,A^{-1}\boldsymbol{ {L}}_{p} (\pos),
{\bf{f}}\cdot{\bf{w}}={\boldsymbol{\lambda}}_c\cdot\boldsymbol{ {L}}_{c} ,%= {\bf{f}}^t \,A^{-1}\boldsymbol{ {L}}_{p} (\pos),
\end{eqnarray}
where ${\bf{w}}=[w_1,\dots,w_n]^T$ is the weight vector 
and 
$\boldsymbol{ {L}}_{c}=[ \mathcal{L}_{\phi(r_1)}(\pos_0),\dots, \mathcal{L}_{\phi(r_n)}(\pos_0), \mathcal{L}_1]^T$. Now, using equation \eqref{eq:lambdac} 
\begin{eqnarray}\label{eq:Vweight2}
%{\bf{f}}^t\,{\bf{w}}={\boldsymbol{\lambda}}^t\,\boldsymbol{ {L}}_{p} (\pos)= {\bf{f}}^t \,A^{-1}\boldsymbol{ {L}}_{p} (\pos),
{\bf{f}}_c\cdot{\bf{w}}_c={\bf{f}}_c \cdot A_c^{-1}\boldsymbol{ {L}}_{c},
\end{eqnarray}
where ${\bf{w}}_c=[{\bf{w}}^T,w_*]^T$  is the weight vector %augmented by a 
whose last component is irrelevant in the calculation . 
From \eqref{eq:Vweight}, we finally obtain 
\begin{eqnarray}\label{eq:Vweight2}
{\bf{w}}_c= A_c^{-1}\boldsymbol{{L}}_c.
\end{eqnarray}
which is similar to \eqref{eq:lambdac} but applied to the differential operator over the RBFs evaluated at the point $\pos_0$.

There are several types of RBFs, and the choice of the optimal one for a given problem is still an open question. %In the called 
Among the infinitely differentiable RBFs, %we  find 
those used most often are the Gaussian $\phi(r)=\exp(-(\epsilon r)^2)$, 
the Inverse Quadratic $\phi(r)=1/(1+(\epsilon r)^2)$ and
the  Inverse Multiquadric $\phi(r)=1/\sqrt{1+(\epsilon r)^2}$. All these
%which 
contain a free parameter $\epsilon$, called shape parameter, that controls the flatness  of the RBF (the smaller the flatter). These  three  RBFs are examples of positive definite RBFs that guarantees that the interpolation matrix  is invertible \cite{Fasshauer07}. 
It is well known that large values of the shape parameter lead to well-conditioned linear systems, but to an inaccurate approximation of the operator. On the other hand, small values of this shape parameter lead to accurate results but make the condition number of the interpolation matrix large, and hence, the interpolation coefficients may diverge \cite{Buhmann03}. 

\subsection{%Local Finite Differences 
	RBF-FD for Laplace-Beltrami operator} 
\label{subsec:LBforphi}

%From \eqref{eq:Vweight2}, to obtain a a 
To obtain an RBF-FD approximation  to the LBO using \eqref{eq:Vweight2}, we start 
by %computing the result of 
applying the operator to an RBF. %first we need to calculate the action of operator   applied to  a RBF.
The LBO can be defined as the %superficial 
surface divergence of the %superficial 
surface gradient,  
\begin{equation}\label{eq:laplacianotangencial}
\Delta_{\Sigma}(\cdot)=\text{div}_{\Sigma}  \left(\boldsymbol{\nabla}_{\Sigma}(\cdot)\right)=
\boldsymbol{\nabla}_{\Sigma}\cdot\boldsymbol{\nabla}_{\Sigma}(\cdot).
\end{equation}
Here  $\boldsymbol{\nabla}_{\Sigma}$ is the surface $nabla$ operator, which is defined as the orthogonal projection of the usual $nabla$ onto the tangent plane to the surface at each point $\pos\in \Sigma$. Thus,
\begin{eqnarray}\label{eq:nablasuperficial}
{\boldsymbol{\nabla}_{\Sigma}=\boldsymbol{\nabla} -\mathbf {n}(\pos)\, \mathbf{n}(\pos)\cdot \boldsymbol{\nabla}}
\end{eqnarray}
where  $\nn(\pos)$ is the  unit normal vector to the surface.

We start %calculating 
by computing the surface gradient %applied to a  
of an RBF ${\boldsymbol{\nabla}_{\Sigma}}\phi$. %\phi(r_i(\pos))$. 
Taking into account that  $\phi(r_i(\node))$ only depends on the distance  $r_i(\node)$ and  applying the chain rule, we obtain
\begin{align}\label{eq:gradienterbff0}
\boldsymbol{\nabla}_{\Sigma}\phi(r_i(\pos))&=\frac{\text{d}\phi(r_i(\pos))}{\text{d}r_i(\pos)}\boldsymbol{\nabla}_{\Sigma}(r_i(\pos))=\nonumber \\
&=\frac{\text{d}\phi(r_i(\pos))}{\text{d}r_i(\pos)}\left(\boldsymbol{\nabla} r_i(\pos) -\left(\nn(\pos)\,\nn(\pos)\cdot\nablav\right) r_i(\pos) \right).
\end{align}
%Due
Using the fact that 
%\begin{eqnarray}\label{eq:aux0}
$\boldsymbol{\nabla} r_i(\pos)=\ds\frac{\bfr_i(\pos)}{r_i(\pos)}, $
%\end{eqnarray}
%the 
equation  \eqref{eq:gradienterbff0} can be written as
\begin{align}
\label{eq:gradienterbf}
\boldsymbol{\nabla}_{\Sigma}\phi(r_i)= 
\left(\bfr_i -(\bfr_i\cdot\nn)\nn\right)\ds\frac{1}{r_i}\ds\frac{\text{d}\phi(r_i)}{\text{d}r_i}\,.
\end{align}
Here, and in all that follows,  we write $r_i = r_i(\pos)$ and $\nn = \nn(\pos)$  to simplify the notation.

To complete our calculation
of the LBO of an RBF %, according to 
using \eqref{eq:laplacianotangencial},  we take the surface divergence of \eqref{eq:gradienterbf}, obtaining
\begin{align}
\label{eq:rawnablasigma}
& \nablat\cdot\left(  \left(\bfr_i -(\bfr_i\cdot\nn)\nn\right)\ds\frac{1}{r_i}\ds\frac{\text{d}\phi(r_i)}{\text{d}r_i}\right) = \nonumber \\
&\ds\frac{1}{r_i}\ds\frac{\text{d}\phi(r_i)}{\text{d}r_i}  \nablat\cdot \left(\bfr_i -(\bfr_i\cdot\nn)\nn\right)+
\left(\bfr_i -(\bfr_i\cdot\nn)\nn\right)\cdot\nablat\left(\ds\frac{1}{r_i}\ds\frac{\text{d}\phi(r_i)}{\text{d}r_i}\right) .
\end{align}
For the  first term  of this equation we have
\begin{align}\label{eq:nablasigma01}
&\nablat\cdot \left(\bfr_i -(\bfr_i\cdot\nn)\nn\right)=\nablat\cdot \bfr_i -(\bfr_i\cdot\nn)\nablat\cdot\nn-\nn\cdot\nablat(\bfr_i\cdot\nn),
\end{align}
where we note that the last term is zero because is the scalar product of two perpendicular vectors. Furthermore,
\begin{align}
\nablat\cdot\bfr_i=\nablav\cdot\bfr_i-(\nn\,\nn\cdot\nabla)\cdot\bfr_i=
\nablav\cdot\bfr_i-\nn\cdot(\nn\cdot\nabla\bfr)
=3-1=2
\end{align} 
and 
\begin{align}
\nablat\cdot\nn=\nablav\cdot\nn-(\nn\,\nn\cdot\nabla)\cdot\nn=
\nablav\cdot\nn-\nn\cdot(\nn\cdot\nabla\nn)
=\nablav\cdot\nn
\end{align}
because  
%$\nn\cdot\nabla\nn=\nabla\lVert \nn \lVert^2/2=0$.
%$(\nn\,\nn\cdot\nabla)\cdot \bfr=\nn\cdot(\nn\cdot\nabla\bfr)=\nn\cdot\nn=1$ and $(\nn\,\nn\cdot\nabla)\cdot \nn=\nn\cdot(\nn\cdot\nabla\nn)=\nn\cdot(\nabla\lVert \nn \lVert^2/2)=0$
$\nn\cdot\nabla\bfr=\nn$ and $\nn\cdot\nabla\nn=\nabla\lVert \nn \lVert^2/2=0$

For the second term of equation (\ref{eq:rawnablasigma}), using \eqref{eq:gradienterbf},  
we have
\begin{align}\label{eq:nablasigma02}
(  \left(\bfr_i -(\bfr_i\cdot\nn)\nn\right)\cdot \nablat\left(\ds\frac{1}{r_i}\ds\frac{\text{d}\phi(r_i)}{\text{d}r_i}\right) =\lVert\bfr_i -(\bfr_i\cdot\nn)\nn\lVert^2\ds\frac{1}{r_i}\frac{\text{d}}{\text{d}r}\left(\ds\frac{1}{r_i}\ds\frac{\text{d}\phi(r_i)}{\text{d}r_i}\right)\, .
\end{align}
Noting that $\lVert\bfr_i -(\bfr_i\cdot\nn)\nn\lVert$ is the norm of the projection of $\bfr_i$ over the tangent plane, and applying the Pythagoras Theorem we have 
$\lVert\bfr_i -(\bfr_i\cdot\nn)\nn\lVert^2=r_i^2 +(\nn\cdot\bfr_i)^2.$ 
Finally, collecting the terms (\ref{eq:nablasigma01})-(\ref{eq:nablasigma02}), we obtain
\begin{align}
\label{eq:laplacianorbf}
\Delta_{\Sigma} \phi(r_i)=
\left(1-  \ds\frac{(\bfr_i \cdot\nn)^2}{r_i^2}-\kappa\ds\frac{\bfr_i\cdot\nn}{r_i} \right)\ds\frac{1}{r_i} \ds\frac{\text{d}\phi(r_i)}{\text{d}r_i}
+
\left(1+  \ds\frac{(\bfr_i\cdot\nn)^2}{r_i^2} \right)\ds\frac{\text{d}^2\phi(r_i)}{\text{d}r_i^2}
\, ,
\end{align}
where $\kappa=\nablav\cdot\nn$ is the curvature term.
%, whose absolute value is equal to the absolute value of the classical mean curvature{\footnote{note that this expression is invariant with respect to the sign of the normal, but the sign of the mean curvature depends on which is considered the exterior normal}} of the surface.

In summary, to compute the value of the LBO applied to an RBF at a certain node
%So that, to compute LBO in a determined  node 
$\pos_l\in{\bf{X}}$, $1\le l \le N$, we consider the stencil centered in $\pos_l$ and size $M$, which is defined as a subset of ${\bf{X}}$ consisting of $\pos_l$ and the $M-1$  nearest neighbor nodes. Without loss of generality, we consider that the center of the stencil is $\pos_l=\pos_1$, and $\pos_k$, $k=2,\cdots, M$ the remaining nodes in the stencil.
Then, from \eqref{eq:weight}
\begin{equation}\label{eq:weightLBO}
\Delta_{\Sigma} f(\pos_1)\approx\ds\sum_{k=1}^M w_{k} f _k\, ,
\end{equation}
where the weights are obtained from \eqref{eq:Vweight2} as 
\begin{equation}\label{eq:weightsstencil}
{\begin{bmatrix}
	w_1 \\
	w_2  \\ 
	\vdots  \\ 
	w_M\\ 
	w_*\
	\end{bmatrix}}=
\begin{bmatrix}
\phi(r_{1}(\pos_{1})) & \phi(r_{1}(\pos_{2})) &  \hdots  & \phi(r_{l}(\pos_{M}))&1 \\ 
\phi(r_{{2}}(\pos_1)) & \phi(r_{{2}}(\pos_{2})) &  \hdots  & \phi(r_{{2}}(\pos_{M}))&1 \\ 
\vdots  &  \vdots  &\ddots  & \vdots &  \vdots  \\ 
\phi(r_{{M}}(\pos_1)) & \phi(r_{{M}}(\pos_{2})) &  \hdots  & \phi(r_{M}(\pos_{M}))&1  \\
1 & 1&  \hdots  & 1&0  \\
\end{bmatrix}^{-1}
\begin{bmatrix}
\Delta_{\Sigma} \phi(r_1(\pos_1))\\ 
\Delta_{\Sigma} \phi(r_{2}(\pos_1)) \\ 
\vdots  \\ 
\Delta_{\Sigma} \phi(r_{M}(\pos_1)) \\
0\\
\end{bmatrix}.
\end{equation}
Here, the elements of the vector in the right hand side, $\Delta_{\Sigma} \phi(r_i(\pos_1))$, are computed using (\ref{eq:laplacianorbf}).

%%%%%%%%%%%%%%%%%%%%%%%%%%%%%%%%%%%%%%%
%%%%%%%%%%%%%%%%%%%%%%%%%%%%%%%%%%%%%%%
%
%
%Ahora lo de la superficie
%
%
%%%%%%%%%%%%%%%%%%%%%%%%%%%%%%%%%%%%%%%
%%%%%%%%%%%%%%%%%%%%%%%%%%%%%%%%%%%%%%%
\subsection{Surface characterization}
\label{subsec:surfchar}
\begin{figure}[t]
	\centering
	\subfigure[]{\includegraphics[width=0.4\linewidth]{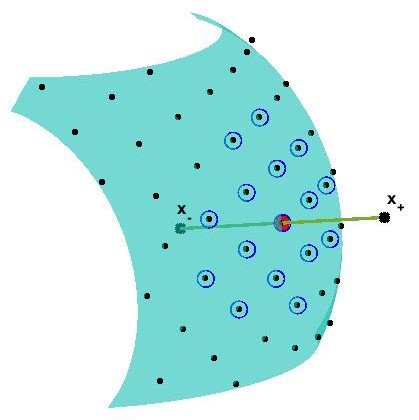}}
	\subfigure[]{\includegraphics[width=0.4\linewidth]{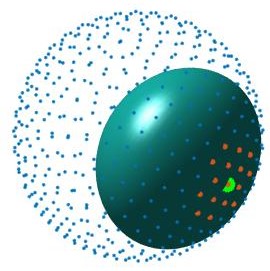}}
	\caption{(a) A closer look to the stencil with the normal vectors, and the on and off points used for interpolation. (b)  Sphere represented by a set of ME blue points. We have also represented a 16 points stencil (in brown) with the central (green) point. We also show the local approximation of the surface $\Sigma_{\Psi}$.} %\textcolor{red}{Poner $\pos_{+}$ y $\pos_{-}$ en la figura.}}
	\label{fig:localsurface}
\end{figure}
When the normal vectors $\nn$ and the curvature $\kappa$ are unknown, it is necessary to compute them
before applying %to apply 
\eqref{eq:laplacianorbf}. %, it is necessary to compute them. 
In this section, we describe a procedure to approximate the values of these quantities.  %normal vectors and curvature. 
The procedure uses a %level set formalism to compute  
% Next, we present a procedure  that use 
level set formalism  and   RBF interpolation to construct  %locally a 
a local approximation to the surface, which %permits us 
is used to compute %these quantities.
their approximate values.% of normal vectors and curvature.  

%In , the authors present a procedure to characterize the surface using a level set formalism with global RBF interpolation. The main advantage of this approach is that one can represent the surfaces on a Cartesian grid without having to parametrize them. This makes the level-set method very flexible when, for example, a surface changes with time. Furthermore, once the level set function $\Psi$ has been found, the unit normal vector and the curvature can be easily computed at any point on the surface.
%\textcolor{red}{ On other hand, both the normal an curvature are local quantities, this is why it is a good idea to adapt this procedure to perform this calculation locally.}

We consider a stencil of size $M$ centered at the node where we want to compute
% that we need to know 
$\nn$ and $\kappa$. %In the same way 
Assume, as before, that $\pos_1$ is the center of the stencil and $\pos_k$, $k=2,\cdots,M$, the other nodes of the stencil.
%First we  find a 
We start by defining a continuously differentiable level-set function $\Psi:\R^3\rightarrow \R$, 
such that its zero level set
\begin{equation}
\label{eq:levelset}
\Sigma_\Psi = \{\node\in \R^3 : \Psi(\node) = 0\},
\end{equation}
approximates locally the surface $\Sigma$ in a neighborhood of $\pos_1$.  %For this, 
To this end, we use an RBF interpolant  \eqref{eq:interpt2} %imposing the  
with the interpolant conditions $\Psi(\node_{k}) =0,   k = 1, \dots , M$.  %In order to avoid the trivial solution  $\Psi(\node) =  0$, we 
We also add two off-surface points $\pos_{+}$ and $\pos_{-}$,  where $\Psi(\pos)$ has non zero values. The resulting interpolant will then be zero at the nodes of the stencil and non zero elsewhere. The off-surface points are placed at either side of the surface along the direction normal to the surface at  $\pos_1$ (see Figure \ref{fig:localsurface}(a)).
Although we assume that the normal vectors are not known at this stage, we can approximate them %roughly 
as 
$$\nn_{app}=(\pos_1-\pos_a)\times(\pos_1-\pos_b), $$ 
%being 
where $\pos_a$ and $\pos_b$ are two  points of the stencil not aligned with $\pos_1$.
% $\pos_{l(1)}$ ,checking that  $\pos_1$, $\pos_a$ and $\pos_b$ are not aligned.

Thus, we model the local surface implicitly as %with a level-set function $\Psi_l$  considering the RBF-interpolant (\ref{eq:interpt}),
\begin{equation}\label{eq:interpolant21}
%\Psi_l(\pos)=\sum_{k=1}^{m} b_i\phi(r_{\bfu(k)}(\pos))+b_{m+1}\phi(\Vert\pos-\pos_+\Vert)    +b_{m+2}\phi(\Vert\pos-\pos_-\Vert)      m{}+ d, \nonumber\\
\Psi_l(\pos)=\sum_{k=1}^{M+2} b_k\phi(r_k(\pos))+ d, \
\end{equation}
%with the constraint 
%\begin{equation}
%\ds\sum_{k=1}^{M+2}b_k=0.\nonumber\ \label{eq:ligadura2}
%\end{equation} 
%In equation \ref{eq:interpolant21}, we have defined 
where $r_{M+1}=\Vert\pos-\pos_+\Vert$ and $r_{M+2}=\Vert\pos-\pos_-\Vert$.
%  to simplify the expression
The interpolation coefficients $b_k$ and the constant $d$ are obtained %from \eqref{eq:interpolacion_c} with    
by enforcing the interpolation conditions (see Figure \ref{fig:localsurface}(b))
\begin{eqnarray}\label{eq:RBFlevelset} 
\Psi_l(\node_{\bfu(k)}) = &0,   \quad k = 1, \dots , M \,\, \qquad \qquad  \quad &\mbox{(on-surface nodes)}\, , \nonumber \\
\Psi_l(\node_+)=  & 1 ;  \Psi_1(\node_-)  = -1 \qquad \qquad  \quad &\mbox{(off-surface nodes)}.
\end{eqnarray}
%as interpolation conditions,  see Figure \ref{fig:localsurface} (b).
and the constraint 
\begin{equation}
\ds\sum_{k=1}^{M+2}b_k=0.\nonumber\ \label{eq:ligadura2}
\end{equation} 

%Next, o
Once the function  $\Psi(\node)$ is known, the unit normal vector to the surface $\Sigma_{\Psi}$  at a point $\node$  is simply given by
\begin{equation}
\label{eq:normal}
\normal = \frac{\boldsymbol{\nabla}\Psi}{\lVert \boldsymbol{\nabla}\Psi\rVert}\, ,
\end{equation}
where the gradient can be computed applying the chain rule
\begin{equation}
\label{eq:gradpsi}
\boldsymbol{\nabla}\Psi = \sum_{k=1}^{M+2} b_k \,\frac{\text{d}\phi(r_k)}{\text{d}r_k}\,
\frac{\bfr_k}{r_{k}}\, .
\end{equation}
%With de normal,  t
Similarly, the curvature  can be obtained from $ \kappa=\boldsymbol{\nabla}\cdot\nn$. Note that
\begin{align} \label{eq:curvatura11}
\boldsymbol{\nabla} \cdot\nn=&
\frac{\lVert\boldsymbol{\nabla} \Psi\rVert\,\boldsymbol{\nabla}\cdot\boldsymbol{\nabla} \Psi-\boldsymbol{\nabla}\Psi\cdot\boldsymbol{\nabla}\lVert\boldsymbol{\nabla} \Psi\rVert }{\rVert\boldsymbol{\nabla} \Psi \lVert^2} %\nonumber  \\
=& \frac{1}{\lVert\boldsymbol{\nabla} \Psi \rVert} \left(  \Delta \Psi-
\nn\cdot\boldsymbol{\nabla} \left(\boldsymbol{\nabla} \Psi \right) \right)\, ,
\end{align}
where we have used \eqref{eq:normal} and %that 
the relation
$\boldsymbol{\nabla}\lVert \boldsymbol{\nabla}\Psi\lVert^2=
\boldsymbol{\nabla}(\boldsymbol{\nabla}\Psi\cdot\boldsymbol{\nabla}\Psi)$. % and \eqref{eq:normal}.
%This expression can be written in a definitive form as 
Then, using \eqref{eq:laplacianorbf} and \eqref{eq:gradienterbf} we obtain
\begin{equation}{\label{eq:curvaturarbf}}
\kappa=
\frac{1}{\lVert \boldsymbol{\nabla} \Psi\lVert}
\ds \sum_{k=1}^{M+2} \ds b_k \left( 
\left(1+\frac{(\bfr_k\cdot\nn)^2}{r_{k}^2} \right)\frac{1}{r_{k}}\ds\frac{\text{d}\phi(r_k)}{\text{d}r_k}
+
\left(1 - \frac{(\bfr_k\cdot\nn)^2}{r_{k}^2} \right) \ds\frac{\text{d}^2\phi(r_k)}{\text{d}r_k^2}
\right) \, ,
\end{equation}
with $\boldsymbol{\nabla} \Psi $ given by \eqref{eq:gradpsi}. %}
Evaluating \eqref{eq:normal} and \eqref{eq:curvaturarbf} at $\pos_1$ we obtain  approximations to the normal vector and the curvature %term in the center of the stencil
at each node, respectively. %Note that, although we do not have control of the direction of  $\nn_{app}$, the equation \eqref{eq:laplacianorbf} do not depends on this.

%
%\textcolor{red}{This formula for the curvature  depends on the normal vectors and on the first and second derivatives of the RBFs at the points that define the surface. 
%Its accuracy only relies on the quality of the interpolation.
%Hence, $\kappa=\kappa(\node)$ depends on the normal vector and on the first and second derivatives of the RBFs at the
%point $\node$ on the interpolated surface. Its accuracy only relies on the quality of the interpolation \eqref{eq:interp3}.}
%
%EDP RBF-FD
%
%\subsection{Solving Surface Equation: Stability }
\subsection{RBF-FD solution of reaction-diffusion equations at surfaces }
\label{stability}
Reaction-diffusion  PDEs on surfaces  %can be represented in the  form
are described by
\begin{equation}\label{eq:reacc_dif_eq}
\frac{\partial }{\partial t} u=R(u,t)+D {\Delta_{\Sigma}} {{u}} \quad\text{on} \quad\Sigma, 
\end{equation}
supplemented with appropriate boundary and initial conditions. Here, 
%$\Sigma$ is the surface on which the equation must be solved, 
$u=u(\pos, t)$ represents the variable of interest,  $R(u,t)$ 
the reaction term, %takes into account the local  variation  laws of $u$, 
$D$  the diffusion coefficient, %which %modulates the surface  diffusion mathematically 
and ${\Delta_{\Sigma}} {{u}} $ the LBO that describes diffusion on the surface $\Sigma$.
%
%multiplies the diffusion term which is described by by the Laplace-Beltrami Operator,   ${\Delta_{\Sigma}} {{u}} $.

%A reactionÐdiffusion system can be solved by using methods of numerical mathematics. There are existing several numerical treatments in research literature.[40][20][41] Also for complex geometries numerical solution methods are proposed.[42][43]

Let  ${\bf{u}}(t)$  be a vector of $N$ components  containing the values of $u(\pos,t)$ at the points     $\pos_i\in\Sigma$,  $i=1,\cdots,N$. % the e
Equation \eqref{eq:reacc_dif_eq} can be spatially discretized as
\begin{equation}\label{eq:reacc_dif_disc}
\frac{d}{dt} {\bf{u}}=R({\bf{u}},t)+D M_{\Delta_{\Sigma}} {\bf{u}}\, ,
\end{equation}
where $M_{\Delta_{\Sigma}} $ is the differentiation matrix %with 
containing the RBF-FD weights for the LBO. %Therefore, the  e
Equation \eqref{eq:reacc_dif_disc} %can be interpreted as 
is a system of $N$ coupled ordinary differential equation and, provided it is stable, it can be advanced in time with a %suitably chosen 
suitable time-integration method.

To %calculate the  
compute the matrix  $M_{\Delta_{\Sigma}} $ we start by selecting the size of the stencil at each node, and use \eqref{eq:normal} and \eqref{eq:curvaturarbf} to obtain $\nn(\pos_i)$ and $\kappa(\pos_i)$, $i=1,\cdots,N$. % at each node. %,  if the normal and curvature are unknown, firstly we compute these quantities for all nodes.  
%Then we select the size of the stencil, and we apply equation \eqref{eq:normal} and \eqref{eq:curvaturarbf} to obtain $\nn(\pos_i)$ and $\kappa(\pos_i)$, $i=1,\cdots,N$.
%%firstly we computed $\nn(\pos_i)$ and $\kappa(\pos_i)$ usung , respectivaly. We choice an appropriate size of the stencil depending on the precision required, the points distribution on the surface, the shape of $\Sigma$.
%Known $\nn(\pos_i)$ and $\kappa(\pos_i)$, we choose the size of the stencil M to apply \eqref{eq:weightLBO}. 
Then,  for each node $\pos_i$, we create  a vector $\bfv$ of $n$ components containing the indexes of the nodes in the stencil %nodes, with 
($\bfv(1)=i$ and $\bfv(k)$,  $k=2,\dots,M$, the indexes of the nearest neighbors nodes to $\pos_i$), and use 
%. Using 
equation \eqref{eq:weightsstencil} to compute %we calcule 
the RBF-FD weights $w_{i,\bfv(k)}$, $k=1,\cdots,M$. Each weight $w_{i,\bfv(k)}$ %constitute the 
is stored as element ($i,\bfv(k)$) of matrix $M_{\Delta_{\Sigma}}$. 
%Repeting this process with each node, $\pos_i$, $i=1,\cdots,N$, we set up LBO matrix. 

A very relevant property that has to be considered is the stability of the time integration method used to solve (\ref{eq:reacc_dif_eq}). In general, a minimum requirement for stability is that all eigenvalues of the differentiation matrix $M_{\Delta_{\Sigma}}$ lie in the left half-plane of the complex plane.  
%Previous to use a  method to temporal integration of \eqref{eq:reacc_dif_disc}, it is necessary to study the stability of this equation. A stability linear study conclude that is $M_{\Delta_{\Sigma}}$ has eigenvalues with real part positive, the equation \eqref{eq:reacc_dif_disc} will not integrable. The m
$M_{\Delta_{\Sigma}} $ is a sparse $N\times N$ matrix, % but 
%thanks to the fact that we are using a local RBF method, it is sparse, with only 
with $M$ entries different from zero per row. Unfortunately, the construction of this matrix does not guarantee that all its eigenvalues lie on the left half-plane and, therefore, stability is not ensured. Fortunately, however, we have found that, if the nodes are suitably chosen and if the internodal distance is small enough, then the eigenvalues do in fact lie in the left half-plane.
%In general, $M_{\Delta_{\Sigma}} $  has not properties such as symmetry or diagonal domminace, hence, we can not know, {\it{a priori}}, if its  eigenvalues are or not real and negative. It is necessary compute these eigenvalues for each case.
Finding the conditions that could guarantee this property is still an open problem. In this context,  the value of the shape parameter and the distribution of the nodes are crucial
\cite{Fuselier13, Shankar15}.
%\textcolor{red}{alguna refeencia}
%The value of the shape parameter is important because, as mentioned above, it has to be small enough to get a good approximation to the operator but not too small to avoid a large condition number of the interpolation matrix $A$ which would result in a wrong set of RBF-FD weights (see for example \cite{Grady}). 
%When possible, the choice of the set of nodes over which we are going to solve the equation is also very important. In a way it is related with the shape parameter, because if the distribution of nodes is not regular, that is if the internodal distance varies a lot, the choice of the shape parameter...

\section{Numerical tests}
In this section, we present  numerical experiments %in order 
whose aim is to test the performance of the proposed procedure to  approximate  the LBO. These experiments are focused (i) on the convergence of the numerical approximation to the LBO applied to a function at a point, (ii) the quality of the characterization of the surface and (iii) the stability properties of the differentiation matrix resulting from the method described in the previous section. All numerical experiments have been perform using $N$ scattered points $\{\pos_i\}_{i=1}^N$ distributed  on the surface of the unit sphere
$\Sph=\{(x,y,z)\in \mathcal{R}^3:x^2+y^2+z^2=1\}$. The used scattered points  are uniformly distributed  %have been selected following the 
Minimal Energy (ME)  nodes, which correspond to the equilibria locations of mutually repelling particles. The exact locations of these ME nodes have been obtained %distribution 
from the repository  \cite{repositorionodos}. In all the experiments we use Gaussians as RBFs.

\subsection{Convergence}

In the first test, we estimate the error of the numerical approximation to $\Delta_{\Sph}f(\pos)$ 
%\begin{equation}
%\label{eq:laplaciano}
%\Delta_{\Sph} f(\pos) = \sum_{i=1}^{M} w_i f(\pos_i)
%\end{equation}
%where the weights $w_i$ were obtained 
using equation \eqref{eq:weightLBO}. % \eqref{eq:Vweight2}.
% for the Laplace-Beltrami operator over the sphere $\Delta_{\mathcal{S}^2}$. 
We apply this operator to the function  $f(x,y,z)=x(1+y(1+z))$, whose exact $\Delta_{\Sph}$ is %can be calculated analytically
\begin{equation}
\label{eq:probanalitico1}
{\Delta} _{\Sph}f=-2x(1+3y(2+z)).
\end{equation}
Note that in the  case  of the unit sphere, if  $\pos\in\Sph$ %   we  have 
then $\nn(\pos)=\pos$ and $\kappa(\pos)=\nabla\cdot\nn(\pos)=2$ (see \cite{Alvarez18} for details).  
%Let $\pos_i$ be a point where we want evaluate the Laplace-Beltrami approximation using  \eqref{eq:laplaciano}, we choose the $M-1$ nearest neighbors to form an $M$-size node stencil. 
The error at location $\pos_1$ is given by
\begin{equation}\label{error}
\text{Error}=\vert {\Delta} _{\Sph}f(\pos_1)-\ds\sum_{i=1}^{M}w_i f(\pos_i)\vert
\end{equation}
where $\pos_i, i=2, \dots M$ are the locations of the $M-1$ nearest neighbors to $\pos_1$.

\begin{figure}[htbp]
	\centering
	\subfigure[]{\includegraphics[width=0.4\linewidth]{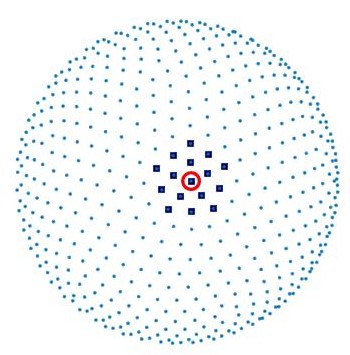}}
	\subfigure[]{\includegraphics[width=0.5\linewidth]{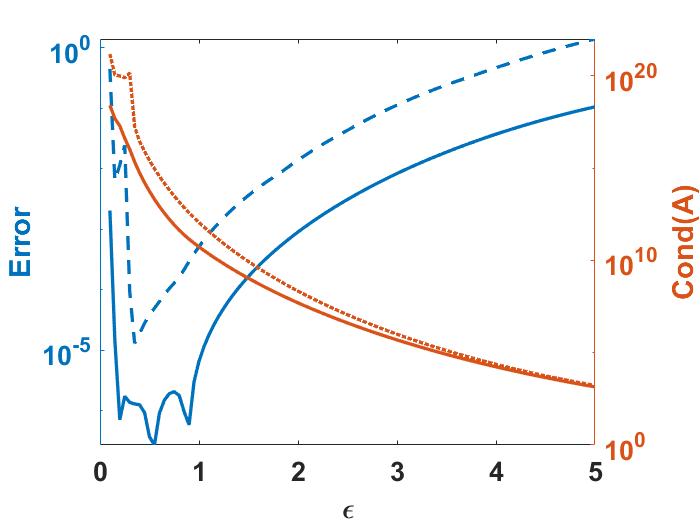}}
	\caption{(a) 1000 ME node distribution over the sphere. The nodes of a  $M=16$ nodes stencil are marked in black, and the center with a red circle.(b) Error and condition number of the interpolation matrix. Continuous lines: data for the point drawn in a) dashed lines:  maximum error and condition number for all nodes. }
	\label{fig:stencil}
\end{figure}

For this experiment we have used $N=1000$ ME nodes which are shown in 
%In 
Figure \ref{fig:stencil}(a). The stencil center is $\pos_1=(1,0,0)$ (marked with a red circle in the Figure), and the stencil size is $M=16$. The stencil nodes are %  and the center ,  we show the points distributions  for  $N=1000$ ME nodes, a stencil of size $16$ nodes,
% marked 
shown with thick blue dots. % in the Figure. %, and the center of stencil , $\pos_1=(1,0,0)$,  marked with a red circle. For this points and with this stencil we calculate the error, 
%\begin{equation}\label{error}
%\text{Error}=\vert {\Delta} _{\Sph}f(\pos_1)-\ds\sum_{i=1}^{16}w_i f(\pos_i)\vert
%\end{equation}

Figure \ref{fig:stencil}(b) shows the resulting error \eqref{error} with a continuous blue line, and the condition number of the interpolation matrix $A_c$ with a continuous brown line, as a function of the shape parameter $\epsilon$.
%and the condition number of the interpolation matrix $A_c$.
%The results are plotted in Figure \ref{fig:stencil}(b) (continuous lines). 
As expected, when the value of $\epsilon$ decreases, the condition number increases and the error decreases until we reach a value of $\epsilon$ where the interpolation matrix becomes ill-conditioned and the error begins to increase sharply. 

If instead of considering a single node and its corresponding stencil we repeat the experiment for all the $N=1000$ nodes and compute the maximum error and maximum condition number, we obtain the results shown with dashed lines in the Figure. %We also plot in  Figure \ref{fig:stencil}(b) (discontinuous lines)  the maximum error and maximum condition number when the center of stencil moves over all points, $\pos_i$, $i=1,\cdots,N.$  And, w
We observe that the behavior is similar to that shown for the case of a single point % point $\pos_1$.  
but the errors and condition numbers are higher.
As we can see,  there is an optimal value of the shape parameter for which the error is minimum. %  the  compromise solution between accuracy and conditioning.

\begin{figure}[htbp]
	\centering
	\subfigure[]{\includegraphics[width=0.45\linewidth]{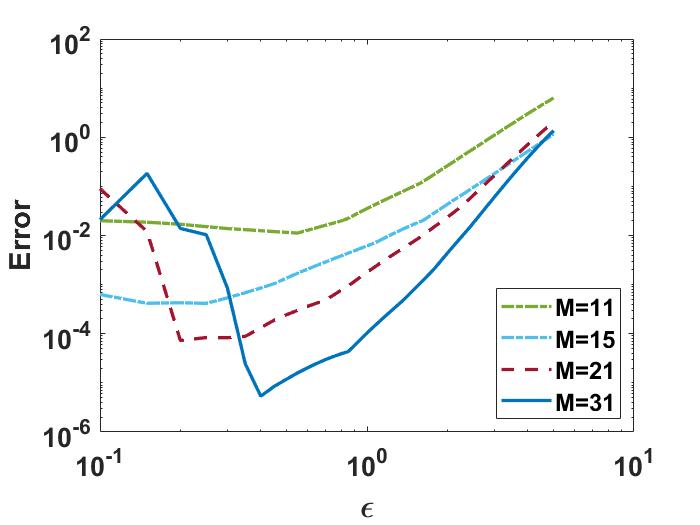}}
	\subfigure[]{\includegraphics[width=0.45\linewidth]{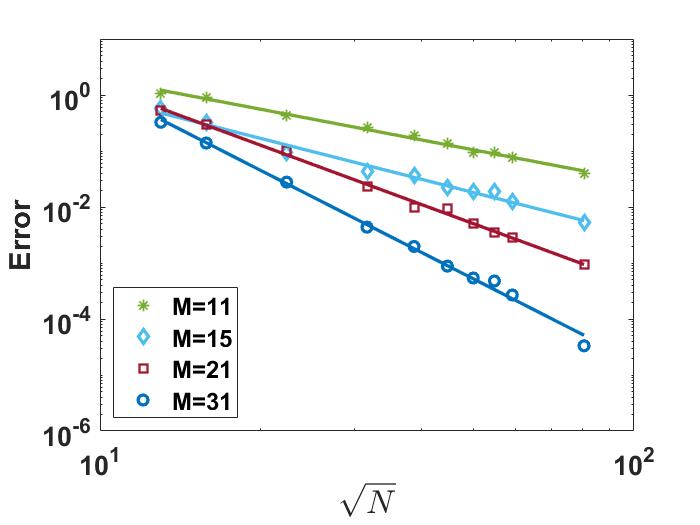} }
	\caption{ (a) Error  for $N=1000$ ME nodes as a function of $\epsilon$ for different values of $M$. (b) Error for $\epsilon=2$ as a function of $\sqrt{N}$ for different values of $M$. %, and we show the error as a function of the number of total points. We have done this for different stencil sizes. We have also matched the results with order one polynomials to check the convergence rate, obtaining 
		Also shown are straight lines with slope $- \mu$ which best fit the data.
		%$\mu=  0.9,1.2,1.7,2.7$  for $M=11,15, 21,31$ stencils size.
	} 
	\label{fig:convergence}
\end{figure}

%In order t
To analyze the order of convergence %$\mu$ (Error $\propto\left(\sqrt{N}\right)^{-\mu} $) 
of the method we have numerically analyzed, for different values $M$ of the stencil size, the dependence of the error on the shape parameter $\epsilon$ and on the internodal distance, which is inversely proportional to $\sqrt{N}$. These results are shown in Figures \ref{fig:convergence}(a) and \ref{fig:convergence}(b).
% and on the %studied the error depending on the density of points and 
%size of the stencils, for a specified value of the shape parameter $\epsilon$. As we can see in 
We use $N=1000$   in Figure \ref{fig:convergence}(a), 
and $\epsilon=2$ in Figure \ref{fig:convergence}(b). %shows the error increasing the size of stencil decreases the  error.  
%The 
Notice that the optimal shape parameter  increases with increasing stencil size. Also
notice that the convergence rate with $\sqrt{N}$ is algebraic, %it seems that the behavior of the convergence is algebraic 
while in the global method the convergence is spectral \cite{Alvarez18}.  To confirm this 
fact we compute the value of the order $\mu$ that best matches the numerical results
%$\mu$ 
(Error $\propto\left(\sqrt{N}\right)^{-\mu} $). %fact, in figure \ref{fig:convergence}(b) we have represented the  error as a function the point density for different stencil sizes. 
We find $\mu=  1.8,2.4,3.4,5.4 $ %$\mu=  0.9,1.2,1.7,2.7$  
for stencils sizes $M=11,15, 21,31$.
As expected, the convergence rate  increases with the stencil size.

%
%(11)   -0.9102    4.8670
%(15) -1.2125    5.4791
%(21)-1.7576    8.4635
%(31)-2.6976   12.6386

\subsection{Surface characterization}
In this Section, we analyze the accuracy in the computation of %results obtained when we calculate 
the normal vector \eqref{eq:normal} and the curvature \eqref{eq:curvaturarbf} using the procedure described in Section \ref{subsec:surfchar}.
For simplicity, we carry out this analysis for the  case  of the unit sphere $\Sph$, for which the exact values of the normal vector and curvature are known. %comparing  with the known exact value. 

\begin{figure}[htbp]
	\centering
	\subfigure[]{\includegraphics[width=0.45\linewidth]{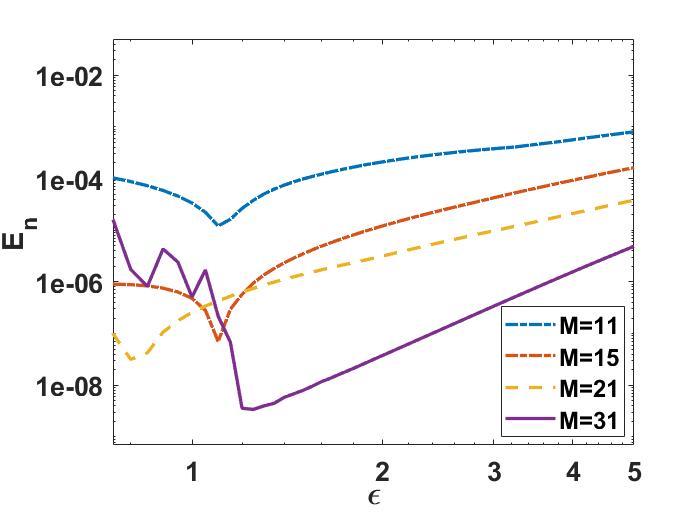}}
	\subfigure[]{\includegraphics[width=0.45\linewidth]{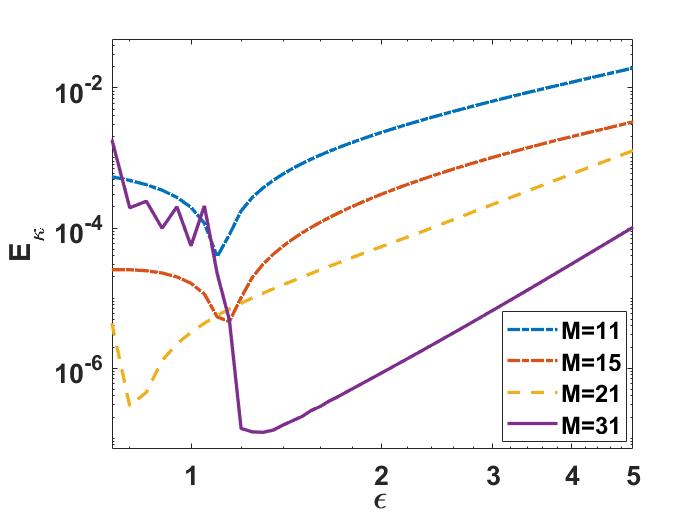}}
	\subfigure[]{\includegraphics[width=0.45\linewidth]{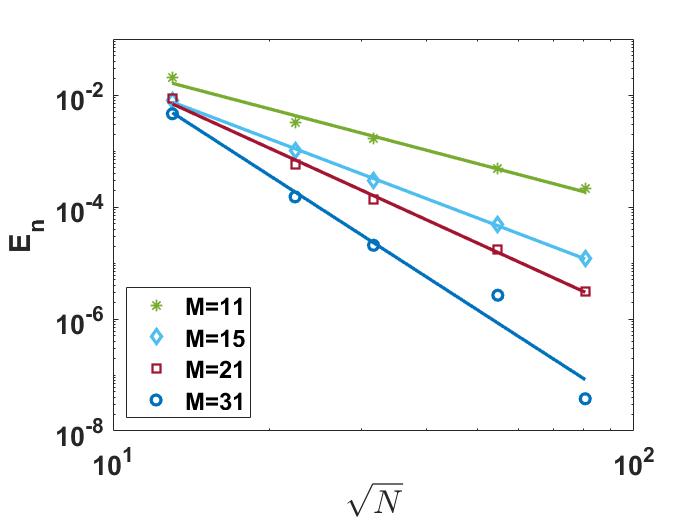} }
	\subfigure[]{\includegraphics[width=0.45\linewidth]{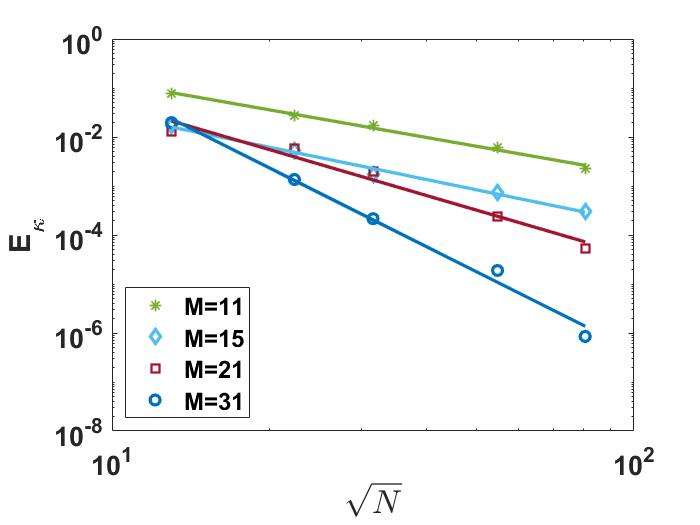} }
	\caption{(a)  Maximum error in the infinity norm of the normal vector as a function of $\epsilon$.
		%	Error en la normal con 
		$N=1000$  ME nodes.  (b)  Maximum error in curvature as a function of $\epsilon$. $N=1000$  ME nodes.
		%Error en la curvatura con $N=1000$  nodos ME.
		(c) Maximum error in the infinity norm of the normal vector as a function of $\sqrt{N}$. $\epsilon =2$. Also shown are straight lines with slope $- \mu$ which best fit the data.
		%Error en la normal con $\epsilon=2$, $\mu=1.2,1.7,2.1,3.0$ para $M=11,15,21,31$, respectivamente.  
		(d) Maximum error in curvature as a function  of $\sqrt{N}$. $\epsilon =2$. Also shown are straight lines with slope $- \mu$ which best fit the data.
		%Error en la curvatura con $\epsilon=2$, $\mu=0.9,1.0,1.5,2.6$ para $M=11,15,21,31$, respectivamente.
	}
	\label{fig:errores_superficie_s2}
\end{figure}

%In f
Figure \ref{fig:errores_superficie_s2} (a) shows the maximum error in the infinity norm of the normal vector, $E_{\nn}=\text{max}\{\Vert \nn(\pos_i)-\pos(i)   \Vert_{\infty}\}$, $ i=1,\cdots, N$, as a function of the shape parameter.
Figure \ref{fig:errores_superficie_s2} (b) shows the maximum error in the computation of the curvature,  $E_{\kappa}=\text{max}\{\vert \kappa(\pos_i)-2   \vert\}$,   $ i=1,\cdots, N$, as a function of the shape parameter.
%$E_{\nn}=\text{max}\{\Vert \nn(\pos_i)-\pos(i)   \Vert_{\infty}\}$ and $E_{\kappa}=\text{max}\{\vert \kappa(\pos_i)-2   \vert\}$,   $ i=1,\cdots, N$, as normal and curvature errors, respectively.
%For these two f
In both Figures we  use $N=1000$ ME nodes.
Again, the error decreases with decreasing $\epsilon$ until the condition number of the interpolation matrix becomes very large and round-off errors deteriorate the accuracy.
% of the interpolation matrix becomes too large to calculate good interpolation weights. %That is why the behavior of these errors is similar to the one of the error in the $\Delta_{\Sph}$ calculation in figure \ref{fig:stencil}. We can also appreciate 
Also observe that the error decreases with increasing stencil size. % improves the error.

%In figures 
Figures \ref{fig:errores_superficie_s2} (c) and (d) show the errors in normal vector and curvature ($E_{\nn}$ and $E_{\kappa}$, respectively) as a function of $\sqrt{N}$.
We observe that, similarly to what happened with the error in the approximation to the LBO, see Figure  \ref{fig:convergence}(b), convergence is algebraic.  
We have computed the value of the order $\mu$ that best matches the numerical results. In the case of the error in the approximation to the normal vector, we have obtained $\mu=2.4,3.4,4.2,6.0$ %$\mu=1.2,1.7,2.1,3.0$ 
for stencil sizes $M=11,15,21,31$, respectively. Lines with these slopes have been plotted in Figure \ref{fig:errores_superficie_s2} (c). In the case of the curvature error, we have obtained $\mu=1.8,2.0,3.0,5.2$ %$\mu=0.9,1.0,1.5,2.6$ 
for the same stencil sizes. Lines with these slopes are also plotted in this figure.
%are plotted the results of the study of convergence rate for both calculation $\nn$ and $\kappa$. As we can appreciate the behavior is similar to showed in figure \ref{fig:convergence}. 
%($N=169, 250, .... 6000$ ME points)
% estos son los datos del ajuste
%Para la normal
%(11)   -1.2255    2.1643
%(15)   -1.7674    4.1705
%(21)   -2.1229    5.921
%(31)   -3.0073   10.0958
%Para la curvatura
%(11)   -0.9368    2.286
%(15)   -1.0960    1.4843
%(21)   -1.5569    4.1338
%(31)   -2.6708    9.950

\subsection{Stability analysis}
\label{subsection_stability}
\begin{figure}[htbp] %\label{fig:estabilidad}
	\centering
	\subfigure[]{\includegraphics[width=0.45\linewidth]{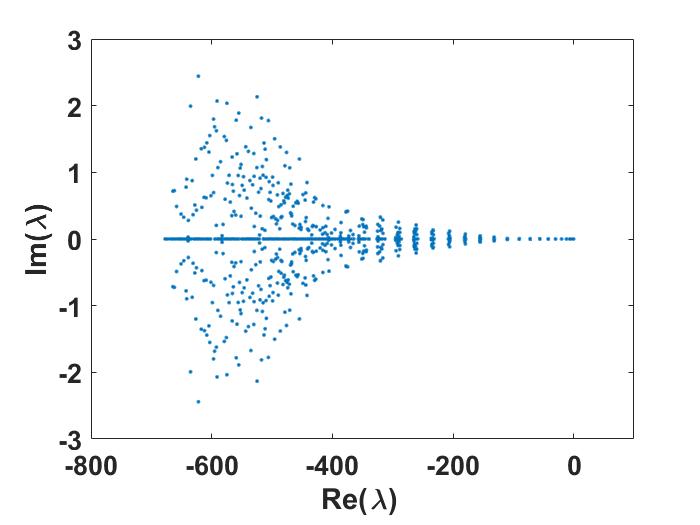}}
	\subfigure[]{\includegraphics[width=0.45\linewidth]{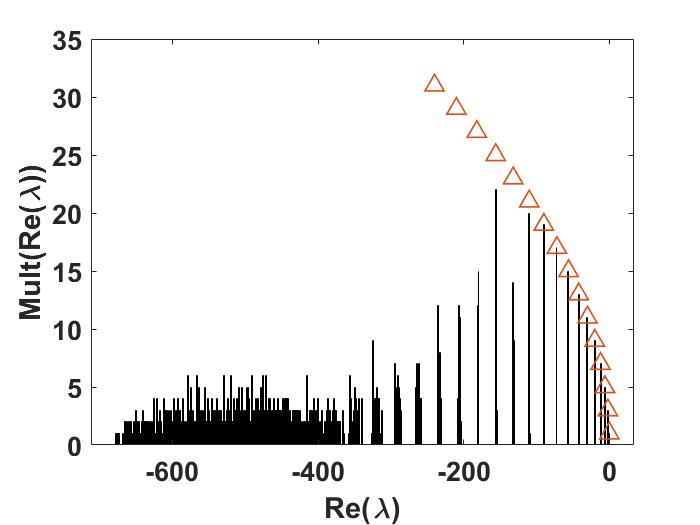}}
	\caption{ (a) Real and imaginary parts of eigenvalues. (b) Histogram of real part of eigenvalues. Triangles represent the exact multiplicity of the eigenvalues given by \eqref{eq:multiplicity}. % multiplicidad teorica de cada uno de ellos. 
	}
	\label{fig:estabilidad}
\end{figure}

As mentioned
in Section \ref{stability}, the eigenvalues of the differentiation matrix  determine the stability of the numerical time integration 
%we have mentioned the importance of the eigenvalues of the matrix $M_{\Delta_{\Sigma}}$ in order to solve numerically 
of reaction-diffusion equations over a surface. Here, we consider the differentiation matrix $M_{\Delta_{\Sph}} $ corresponding to the unit sphere for the case of $N=1000$ ME nodes, $M=31$ stencil size, and $\epsilon=2$ shape parameter.
%Next, we are going to study the eigenvalues of this matrix. The sparsity  of our LBO matrix is showed  in f
%Figure  \ref{fig:estabilidad} (a) shows the sparse pattern of the differentiation matrix where are plotted 
%by plotting the non-zeros entries of $M_{\Delta_{\Sph}} $. % for $N=1000$ ME nodes  and $M=31$as  stencil size.
The exact eigenvalues of the Laplace-Beltrami operator %$\Delta_{\Sph}$ 
are real and non positive; see Chapter 3 of  \cite{Shubin01}. Their exact values are $\lambda_k=-k(k+1)$, $k=0,1,2,\dots$ with multiplicity  %is given by
\begin{equation}
\label{eq:multiplicity}
\textrm{{\it Mult}}(\lambda_k)=\binom{2+k}{2}-\binom{k}{2}.
\end{equation}
Figure \ref{fig:estabilidad} (a) shows the 
eigenvalues of the $M_{\Delta_{\Sph}} $ matrix. %calculated for $N=1000$ ME nodes, $\epsilon=2$,  and $M= 31$ are represented in figure \ref{fig:estabilidad} (b). 
All of them have negative real part, and for small values of $k$ they are very close to the exact ones ($\lambda_k =2, 6, 12, 20, \dots$ for $k =1, 2, 3, 4, \dots$). %are similar to the exact ones as we can see in f
This can be better observed in Figure \ref{fig:estabilidad} (b) which shows the histogram of the real part of eigenvalues (the imaginary part is approximately zero).Since the height of the histogram represents the number of eigenvalues of a certain value, that height represents the multiplicity of the eigenvalue.  The exact multiplicity of the eigenvalues given by \eqref{eq:multiplicity} are
shown with triangles. Notice that, for small values of $k$, the exact multiplicities of the eigenvalues of the Laplace-Beltrami operator are in very good agreement with the multiplicities of the eigenvalues of the differentiation matrix $M_{\Delta_{\Sph}} $.
%ere we have represented the histogram of real part of the eigenvalues along with it theoretical multiplicity, which is also in agreement for small values of $k$.
%
%\begin{figure}[htbp] %\label{fig:estabilidad}
%	\centering
%	\subfigure[]{\includegraphics[width=0.45\linewidth]{figuras/figura20a}}
%	\subfigure[]{\includegraphics[width=0.45\linewidth]{figuras/figura21a}}
%	\subfigure[]{\includegraphics[width=0.45\linewidth]{figuras/figura22a} }
%	\subfigure[]{\includegraphics[width=0.45\linewidth]{figuras/figura21a}}
%  \caption{ (a) Sparsity pattern of matrix $M_{\Delta_{\Sigma}}$. (b) Real and imaginary parts of eigenvalues. (c) Histogram of real part of eigenvalues. Triangles represent the exact multiplicity of the eigenvalues given by \eqref{eq:multiplicity} % multiplicidad teorica de cada uno de ellos. 
%  d) \textcolor{blue}{ya  se vera}}
%  \label{fig:estabilidad}
%\end{figure}

%\textcolor{red}{La Figura \ref{fig:estabilidad} (a) no est� referenciada en el texto. �Quitarla?}

\section{Numerical results}

In this section we apply the method that we propose for the numerical approximation of the LBO %$\Delta_{\Sigma}$ 
in order to compute the solution of two examples of reaction-diffusion problems on surfaces. % where the reaction term is linear and the diffusive term is modeled by $\Delta_{\Sigma}$.

%  En esta seccion  utilizaremos nuestro procedimiento para obtener el  LBO para resolver EPDs sobre superficies.  Los ejemplos eligidos son dos ecuaciones del tipo  reaccion-difusion,  ecuaciones en la que el termino  de reaccion es un termino lineal (\textcolor{red}{por lo que no es muy chuli utilizar metodos implicitos, de esto habra que hablar antes})  y en las que el termino difusivo esta modelado por LBO sobre la superficie en cuestion.  

\subsection{First example:Turing Patterns}

As a first example, we consider the reaction-diffusion system  proposed  by Alan Turing   as a prototype model for pattern formation in nature \cite{Turing52}.  These patterns arise %arising 
% with the evolution with time starting 
from a homogeneous and uniform initial state.
Turing's model can generate a variety of spatial patterns %which can arise from 
which originate from a wide variety of phenomena: morphogenesis in biology \cite{Maini12}, ecological invasion \cite{Holmes94} or tumor growth \cite{Chaplain01}.

%It has also been argued that reactionÐdiffusion processes are an essential basis for processes connected to morphogenesis in biology[21] and may even be related to animal coats and skin pigmentation.[22][23] Other applications of reactionÐdiffusion equations include ecological invasions,[24] spread of epidemics,[25] tumor growth[26][27][28] and wound healing.[29] Another reason for the interest in reactionÐdiffusion systems is that although they are nonlinear partial differential equations, there are often possibilities for an analytical treatment.[8][9][30][31][32][20]and mathematical studies have revealed the kinds of interactions required for each., giving this model the potential for application as an experimental working hypothesis in a wide variety of morphological phenomena. 
%This foundational paper describes the way in which patterns in nature such as stripes and spots can arise naturally out of a homogeneous, uniform state.

%The equations of 
Turing's equations model  the interaction of an activator $u(\pos,t)$ and an inhibitor $v(\pos,t)$. This interaction is described by the following two differential  equations:
\begin{align}\label{eq:turing1}
\frac{\partial u}{\partial t} &=\alpha u (1-\tau_1 u^2)+ v(1-\tau_2 u)+D_{u}\Delta_{{\Sigma}} u\\
\label{eq:turing2}
\frac{\partial v}{\partial t}&=\beta u \left(1+\frac{\alpha\tau_1 u v}{\beta}\tau_1 u^2\right)+ u(\gamma-\tau_2 v )+D_{v}\Delta_{{\Sigma}} v.
\end{align}

Depending on the choice of parameters different patterns can arise naturally out of a homogeneous, uniform state. Both variables, $u$ and $v$, can lead to instabilities which evolve to different pattern formations, such as stripes or spots. %In table \ref{tab:turing} we show the values of the parameters used in our numerical experiments, which coincide with those used in  \cite{Shankar15}. 

The  surface chosen to solve  equations \eqref{eq:turing1}  and \eqref{eq:turing2} is the Schwarz Primitive  Minimal Surface which can be approximated by the implicit equation \cite{websurface}
\begin{equation}\label{eq:SchwarzP}
\cos(2\pi x)+\cos(2\pi y)+\cos(2\pi z)=0.
\end{equation}
We consider periodic boundary conditions in the three coordinate axes, that is, $ u (x = -1, y, z) = u (x = 1, y, z) $ and $ v (x = -1, y, z) = v (x = 1, y, z) $  and the same for the Y- and Z-axes.

\begin{figure}[htbp]\label{fig:turing}
	\centering
	\subfigure[]{\includegraphics[width=0.4\linewidth]{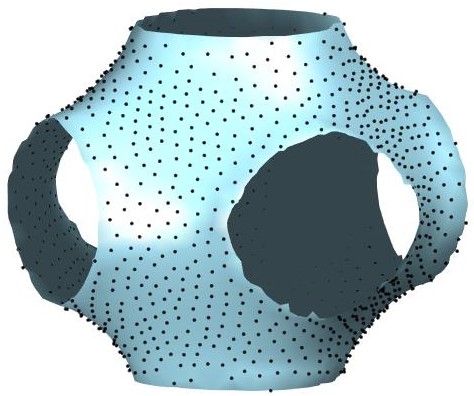}}
	\subfigure[]{\includegraphics[height=1.6in,width=0.4\linewidth]{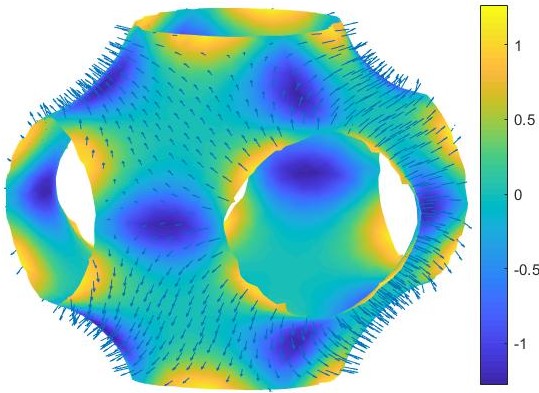}}
	\subfigure[]{\includegraphics[width=0.4\linewidth]{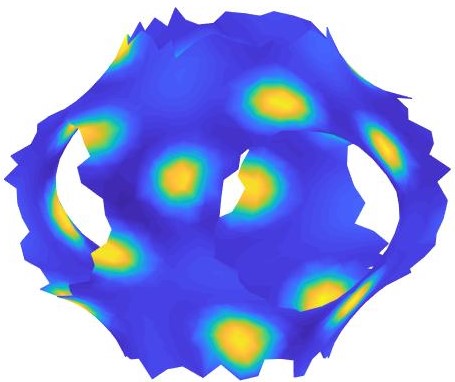} }
	\subfigure[]{\includegraphics[width=0.4\linewidth]{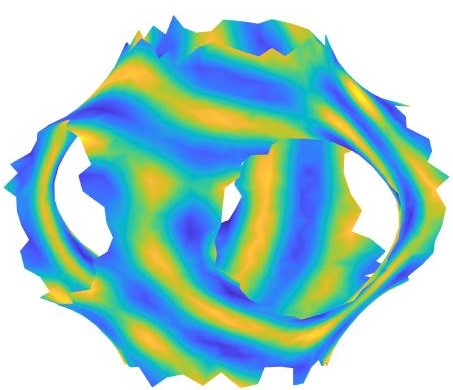}}
	\caption{(a) Nodes and surface. (b) Curvature and normal vectors.(c) Stationary state for variable $u$ resulting from the parameters in the first row of Table \ref{tab:turing}. (d) 
		Stationary state for variable $u$ resulting from the parameters in the second row of Table \ref{tab:turing}.
		% Comentario Manolo\textcolor{red}{que representan los colores de la figura a? no convendr�a poner una colorbar en la figura b para saber los valores de la curvatura?}	
	}
\end{figure}

Figure \ref{fig:turing} (a) shows the surface  \eqref{eq:SchwarzP} and the nodes used for the computation. These nodes have been obtained by projecting radially over Schwarz's surface a set of $N=1800$ ME nodes on the unit sphere obtained from the repository \cite{repositorionodos}. With this set of nodes, we calculate the normal vectors $\nn(\pos_i)$  and the curvature $\kappa(\pos_i)$ on the surface using equations \eqref{eq:normal} and \eqref{eq:curvaturarbf}, respectively. The resulting values  are shown in Figure \ref{fig:turing} (b). %For these calculations we have used 
These values have been obtained using a stencil size $M=51$ and a shape parameter $\epsilon= 3$. Since the analytic formula of the surface is known \eqref{eq:SchwarzP} %As we have an analytic expression for the surface, 
we can calculate the exact values of $\nn(\pos_i)$  and  $\kappa(\pos_i)$ and, thus, the error. We find that these errors are %which has been 
smaller than 0.1$\%$ both for the normal vectors and  for the curvature.

%La distribucion de puntos con lo que hemos resuelto nuestro ejemplo la hemos obtenido proyectando radialmente sobre la superficie \eqref{eq:sup}  los  el conjunto de puntos ME  sobre la esfera unidad obtenidos de \cite{repositoriodenodos}.
% Con dicho conjunto de puntos, en primer lugar caracterizamos la superficie, estimando tanto sus vectores normales,$\nn{\pos_i}$  como su curvatura, $\kappa(\pos_i)$ utlizando las expresiones \eqref{eq:normal} y \eqref{eq:curvatura}. 
% En la figure \ref{fig:turing}(a) hemos representado tanto la curvatura como los algunos vectores obtenido cuando utilizamos un stenicil de tamano XX con un valor parametro de forma $\epsilon=XX$. 
% Con estos valores, dado que esta superficie si se pueden calcual analiticamente estas cosas, los errore decir que los errores  comentido han sido XX $\%$ y XX $\%$, en el computo de las normales y curvatura respectivamente.{\textcolor{red}{esto vere si se puede hacer}}.  

%Once we have calculated 
The approximate values of $\nn(\pos_i)$ and  $\kappa(\pos_i)$ are then used %, we can 
to compute % the approximation to t
the differentiation matrix of the Laplace-Beltrami operator  $M_{\Delta_{\Sigma}}$. % for this surface. In order to do so, 
In this computation we have used  $M=31$ stencils and $\epsilon=2$. 
We have also checked that the eigenvalues of this matrix do not have positive real part.
%Once we have this matrix, we check that it does not have eigenvalues with a positive real part, which ensures stability when we solve the equations \eqref{eq:turing1} and \eqref{eq:turing2} .

For temporal integrations we use  function {\tt ode45} from MatLab, a standard solver for non-stiff ordinary differential equations. %, the function $ode45$ from MatLab. 
This function implements a Runge-Kutta method with a variable time step for efficient computation. 
\begin{table} %[!hbt]
	\begin{center}
		\begin{tabular}{|c|c|c|c|c|c|c|c|}
			\hline
			pattern & $D_u$ & $D_v$ & $\alpha$ &$ \beta$ &$\gamma$& $ \tau_1 $& $\tau_2$\\ 
			\hline
			spots &$ 2.32\,10^{-3}$ & $ 4.5\, 10^{-3}$ &0.899  & -0.91 &-0.899 & 0.02 &0.2 \\  
			\hline 
			stripes &$ 1.08\,10^{-3} $&  $2.1\, 10^{-3}$ &0.899  & -0.91 &-0.899 &3.5 &0  \\  
			\hline
		\end{tabular} 
	\end{center}
	\caption{Parameter values used to solve Turing's model \eqref{eq:turing1}-\eqref{eq:turing2}. %The first set produces a spot pattern, while the second gives a stripes pattern. 
	}
	\label{tab:turing}
\end{table}
%We have solved Turing's problem with the two sets of parameters shown in Table \ref{tab:turing}. These sets have been taken from \cite{Shankar15}. % we show the values of the parameters used in our numerical experiments, which coincide with those used in  \cite{Shankar15}. 
%The first set of parameters produces a spots pattern, while the second results in a stripes pattern. 
The initial condition in both cases is a perturbation of the activator $u$. 
%(\textcolor{red}{deber�amos ser mas especificos}). 
The interaction between the activator and the inhibitor $v$ results in a transient state that, for long times, evolves to %ends in 
a stationary state with %on which appear 
different characteristic patterns which depend on the %election of the equations 
values of the parameters. %In 
Figures \ref{fig:turing} (c) and (d)  show the stationary states corresponding to the two sets of parameters shown in Table \ref{tab:turing},
which have been taken from \cite{Shankar15}. These two stationary states correspond to two  classic patterns of Turing's model: the spots and the stripes pattern, respectively. These stationary states are in agreement with those obtained previously for the same parameters \cite{Shankar15}.

\subsection{Second example: Schaeffer's model}
In the second numerical experiment we solve % apply the proposed method to a 
the bioelectric cardiac source model proposed by Mitchell and Schaeffer \cite{Mitchell02}.
This  is a simple model %, but it is  able to reproduce  
which is capable of reproducing the  main electrophysiological  properties  of cardiac tissue, such as restitution properties or spatial variations of the action potential duration.
The model consists of two differential equations
\begin{equation}\label{eq:propagacion}
\frac{\partial v}{\partial t} =\sigma\Delta_{\Sigma} v+ J_{in}(v,h)+J_{out}(v)+J_{stim}(t)\,,
\end{equation}
and
\begin{equation}\label{eq:chifer2}
\frac{\partial h}{\partial t} =\left\{ \begin{array}{ll}
\displaystyle\frac{1-h}{\tau_{open}},& v < v_{crit} \\
\\
\displaystyle\frac{-h}{\tau_{close}},& v > v_{crit}. \\
\end{array} \right.
\end{equation}
where $v=v(\pos,t)$      and  $h=h(\pos,t)$ are   the transmembrane voltage and the inactivation gate variable, respectively. %Both equations are written in dimensionless form, and the variables are scaled to vary between 0 and 1.
% \textcolor{red}{variables are written in dimensionless form and  scaled to vary between 0 and 1.}
In  \eqref{eq:propagacion}, the terms  $J_{in}$ and $J_{out}$ represent the inward and outward currents of the cells of the membrane
\begin{equation}\label{eq:JinJout}
J_{in}(v,h) = \frac{h(1-v)v^2}{\tau_{in}} 
\quad
\text{and}\quad 
J_{out}(v) = -\frac{v}{\tau_{out}},
\end{equation} 
and $J_{stim}$ represents the  initial stimulus. The diffusive  term $%i_m=
\sigma\Delta_{\Sigma} v$ models the transmembrane current flowing through the cardiac membrane. 
The six parameters %in  \eqref{eq:propagacion}-\eqref{eq:chifer2}  
$\sigma, \tau_{in}, \tau_{out}, \tau_{open}, \tau_{close}$, and $v_{crit}$ govern the behavior of the membrane and, depending on their values, %they take, 
they can not only model healthy tissue, but also tissue with some kind of pathology. We refer the interested reader to \cite{Mitchell02,Alvarez12} for more details about this model.
\begin{table} %[!hbt]
	\begin{center}
		\begin{tabular}{|c|l| c|c|l|c|clc|c||}
			\hline
			${\sigma}$& $\tau_{open}$ & $\tau_{close}$ & $\tau_{in}$ &$ \tau_{out}$ & $ v_{crit} $\\ 
			\hline
			$10^{-3} \, cm^2/ms$ & $130\, ms$ & $150 \, ms$  & $0.2\, ms$&  $10\, ms$& 0.13  \\  
			\hline 
		\end{tabular} 
	\end{center}
	\caption{%Valores de los parametros  utilizados para resolver las ecuaciones 
		Parameters used for the solution of Mitchell and Schaeffer's model \eqref{eq:propagacion}-\eqref{eq:chifer2} . %Dichos parametros se  corresponden con un tejido cardiaco sano.
	}
	\label{tab:chifer}
\end{table}

\begin{figure}[htbp]\label{fig:chifer}
	\centering
	\subfigure[]{\includegraphics[width=0.4\linewidth]{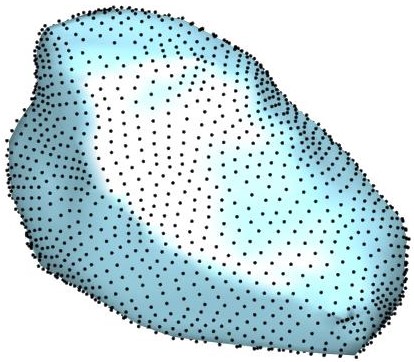}}
	\subfigure[]{\includegraphics[width=0.4\linewidth]{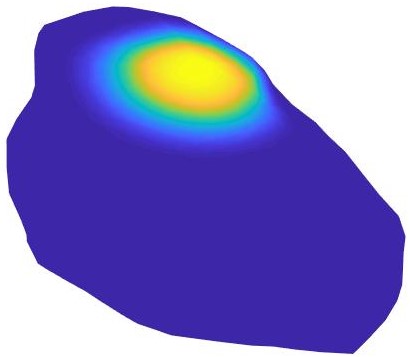}}
	\subfigure[]{\includegraphics[width=0.4\linewidth]{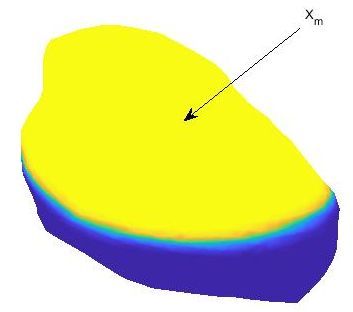} }
	\subfigure[]{\includegraphics[width=0.5\linewidth]{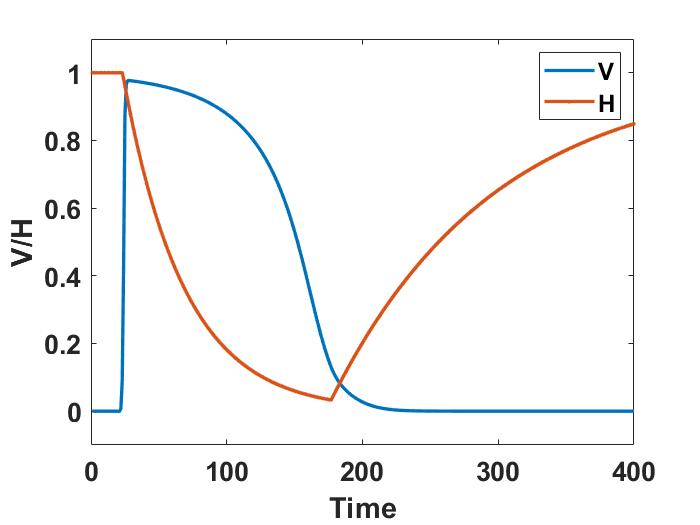}}
	\caption{(a) Nodes and surface. (b) Initial condition for variable $V$.
		% Comentario Manolo \textcolor{red}{que variable se representa?}
		(c) Solution at time $50$ ms for variable $V$.(d) Full heartbeat cycle at point $\pos_m.$}
\end{figure}

We have applied our proposed procedure to the solution of equations \eqref{eq:propagacion}-\eqref{eq:JinJout} with the parameters shown in Table \ref{tab:chifer}. These equations are solved on the epicardium, which is the outermost surface membrane of the heart and which is shown in Figure \ref{fig:chifer}(a). %In the figure \ref{fig:chifer}(a)  are represented 
We also show the $ N =2014 $ points used for the numerical solution of the problem. These points have been obtained from a computerized tomography (CT) of a 
%$ 43 $ years old 
real patient \cite{Chavez15}. Therefore, it is a realistic model in which we do not have information neither on the normal vectors nor on the curvature of the surface.
%Considering
We consider
$$v(\pos,t=0)=0\,  \quad \text{and} \quad h(\pos,t=0) = 1 \, \quad \forall\,\pos\in \Sigma$$  as initial conditions, %an 
and we apply a current stimulus  %\textcolor{red}{(habria que ser mas especificos)}
%is applied, (see Figure \ref{fig:chifer}(b)).
\begin{equation}
\label{eq:stim}
J_{stim}(\pos, t)=H(t_{stim}-t)e^{\frac{(\pos-\pos_s)^2}{\delta^2}}\, ,\quad t\ge 0 \, ,
\end{equation}
shown in Figure \ref{fig:chifer}(b). Here,  $H(t)$ denotes the Heavisade function, $t_{stim}$ the time when the stimulus ends, $\pos_s$ the position where the stimulus is applied, and $\delta$ its spatial width. 

We use  the same method that in the previous example for temporal integration.  The solution shows the propagation of the electric excitation along the membrane. For instance,  Figure \ref{fig:chifer} (c) %, where we have represented 
shows the transmembrane current  50 ms after the stimulus ends. The tissue goes from a  resting to an excited state. %(refractory). 
Finally,  the membrane returns to the  resting state awaiting for the next stimulus. This behavior can be observed in Figure \ref{fig:chifer} (d), where we show %have represented 
the evolution with time of the transmembrane voltage $ v(\pos_m,t)$ and the gate variable $ h (\pos_m,t) $ at the point $\pos_m$ in
Figure \ref{fig:chifer} (c). As we can see,
the cardiac tissue experiences the different stages of a heartbeat corresponding to those shown in \cite{Mitchell02}. %showing the success of the 
%Thus, we show that the proposed methodology is adequate for %proposed here for 
%approximating the solution of this type of PDEs.

\section{Conclusions}
In this article we have introduced a new RBF-FD based method to calculate the solution of reaction-diffusion equations over surfaces. The diffusive part of this kind of equations is usually modeled by the LBO, for which we calculate a numerical approximation, $M_{\Delta_{\Sigma}}$, using an RBF-FD approach. This approach is based on an explicit formula for the LBO applied to an RBF,
which is exact; see (\ref{eq:laplacianorbf}). This formula involves the normal vectors and the curvature to the surface. These are known if
an explicit expression for the surface is available. If the surface is defined by a set of scattered nodes which is the usual output of, for example 3D scanning processes, we propose a level set formalism and an RBF interpolation to construct  a local approximation to the surface, from which we estimate these quantities.

We have study numerically the convergence and stability properties of our approach, and  we have shown the performance of the proposed methodology solving the Turing model for natural pattern formation over a Schwarz Primitive Minimal Surface, and the Schaeffer's model for electrical cardiac tissue behavior. %over a point cloud generated by a computerized tomography of a real patient.

\section{Acknowledgments}
This work has been supported by Spanish MICINN Grant FIS2016-77892-R.


\begin{thebibliography}{10}

\bibitem{Alvarez12}
D.  Alvarez, F.  Alonso-Atienza, J.L.  Rojo-Alvarez, A.  Garcia-Alberola, M. Moscoso,
Shape reconstruction of cardiac ischemia from non-contact intracardiac recordings: a model study,
Math. Comput. Model. 55 (2012) 1770--1781 .

\bibitem{Alvarez18}
D. Alvarez, P. Gonzalez-Rodriguez, and M. Moscoso, A Closed-Form Formula for the RBF-Based Approximation of the Laplace-Beltrami Operator, J. Sci. Comput. (2018) 77: 1115. %https://doi.org/10.1007/s10915-018-0739-1

\bibitem{Barreira11}
 R. Barreira, C.M. Elliott, and A. Madzvamuse, 
 The surface finite element method for pattern formation on evolving biological surfaces, 
J. of Math.Biology. 63 (2011) 1095--1119.

\bibitem{Bertalmio01}
M. Bertalmio, L.-T. Cheng, S. Osher, and G. Sapiro, 
Variational problems and partial differential equations on implicit surfaces, 
J. Comput. Phys. 174 (2001)  759--780.


\bibitem{Buhmann03}
M.D. Buhmann,
Radial Basis Functions,
Cambridge University Press (2003)

\bibitem{Carr01} % esta mejorarla
J. C. Carr, R. K. Beatson, J. B. Cherrie, T. J. Mitchell, W. R. Fright, B. C. McCallum, and T. R. Evans. 2001. Reconstruction and representation of 3D objects with radial basis functions. In Proceedings of the 28th annual conference on Computer graphics and interactive techniques (SIGGRAPH '01). ACM, New York, NY, USA, 67-76. DOI: https://doi.org/10.1145/383259.383266


\bibitem{Chaplain01}
M.A. Chaplain,M. Ganesh and I.G. Graham,  Spatio-temporal pattern formation on spherical surfaces: numerical simulation and application to solid tumour growth, J. Math. Biol. (2001)   42(5):387-423.

\bibitem{Chavez15}
C. E. Chavez,N. Zemzemi, Y. Coudi�re, F.  Alonso-Atienza and D.Alvarez, Inverse Problem of Electrocardiography: estimating the location of cardiac isquemia in a 3D geometry. Functional Imaging and modelling of the heart (FIMH2015), Maastricht,
Netherlands. 9126, Springer International Publishing, 2015,

\bibitem{Dziuk07}
G. Dziuk and C. Elliott, 
Surface finite elements for parabolic equations, 
J. Comp. Math. 25 (2007) 385--407.

\bibitem{Fasshauer07}
G.E. Fasshauer, Meshfree Approximation Methods with Matlab, World Scientific Publishers, 2007.

\bibitem{Floater05}
M. Floater and K. Hormann, 
Surface parameterization: a tutorial and survey, 
Advances in Multiresolution for Geometric Modelling, Springer, 2005.


\bibitem{Fonberg15} 
B. Fonberg and N. Flyer, A primer on radial basis functions with applications to geosciences, CBMS-NSF Regional Conference Series in Applied Mathematics (2015). %ISBN 978-1-611974-02-7.

\bibitem{fornberg2011}
B. Fornberg and E. Lehto, Stabilization of RBF-generated finite difference methods for convective PDEs, Journal of Computational Physics 230 (2011) 2270-2285.

\bibitem{Flyer11}
 N. Flyer and B. Fornberg, 
Radial basis functions: Developments and applications to planetary scale flows, 
Computers and Fluids, 26 (2011) 23--32.

\bibitem{Flyer09}
N. Flyer and G. B. Wright, A radial basis function method for the shallow water equations on a sphere, Proc. Roy. Soc. A 465 (2009) 1949--1976.

\bibitem{Fuselier13}
E. J. Fuselier and G. B. Wright, A high order kernel method for diffusion and reaction-diffusion equations on surfaces, Journal of Scientific Computing (2013) 1-31.

\bibitem{Greer06}
 J. B. Greer, 
 An improvement of a recent Eulerian method for solving PDEs on general geometries,
J. Sci. Comput. 29 (2006)  321--352.


\bibitem{Gu04}
 X. Gu, Y. Wang, T. Chan, P. Thompson, and S. Yau, Genus zero surface conformal mapping
and its application to brain surface mapping, Medical Imaging, IEEE Transactions on, 23
(2004), pp. 949--958.


\bibitem{Holmes94}
 E. E. Holmes, M. A. Lewis, J. E. Banks and R. R. Veit,  Partial Differential Equations in Ecology: Spatial Interactions and Population Dynamics, Ecology (1994) 75(1), 17-29.
 % doi:10.2307/1939378


\bibitem{Holst01}
M. Holst, 
Adaptive numerical treatment of elliptic systems on manifolds, 
Advances in Computational Mathematics  2001 139--191.


\bibitem{Lehto17}
E. Lehto, V. Shankar, and G. B. Wright, 
A Radial Basis Function (RBF) Compact Finite Difference (FD) Scheme for Reaction-Diffusion Equations on Surfaces, 
%J. Sci. Comput. 39 (2017) 2129--2151.
SIAM J. Sci. Comput. (2017) 39 (5)  A2129-A2151. 

\bibitem{Macdonald09}
C. B. Macdonald and S. J. Ruuth, 
The Implicit Closest Point Method for 
the Numerical Solution of Partial Differential Equations on Surfaces,
SIAM J. Sci. Comput. 31 (2009) 4330--4350.


\bibitem{Maini12}
 P. K. Maini, T. E. Woolley, , R. E. Baker,  E. A. Gaffney, and S. S. Lee,   Turing's model for biological pattern formation and the robustness problem, Interface focus (2012) 2(4), 487-496. %doi:10.1098/rsfs.2011.0113


\bibitem{repositorionodos}
https://web.maths.unsw.edu.au/~rsw/Sphere/Energy/index.html %#ME



\bibitem{Mitchell02}
C.C.Mitchell, D. G. Shaeffer, A Two-Current Model for the Dynamics of Cardiac Membrane, B. Math. Biol. 65 (2003), 767-793.

\bibitem{Piret2012}
C. Piret , The orthogonal gradients method: A radial basis functions method for solving partial differential equations on arbitrary surfaces, J. Comput. Phys. (2012) 231, 4662-4675

\bibitem{websurface}
http://www.msri.org/publications/sgp/jim/geom/level/library/triper/index.html



\bibitem{Sarra09}
S.A. Sarra and E.J. Kansa
Multiquadric radial basis function approximation methods for the numerical solution of partial differential equations,
Adv. Comput. Mech., 2 (2009)

\bibitem{Shankar2013}
V. Shankar, G. B. Wright, A. L. Fogelson and R. M. Kirby, A study of different modeling choices for simulating platelets within the immersed boundary method, Appl. Numer. Math. (2013) 63, 58?77

\bibitem{Shankar15}
V. Shankar, G. B. Wright, R. M. Kirby and A. L. Fogelson
 A radial basis function (RBF)- finite difference (FD) method for Diffusion and reaction Diffusion equations on surfaces, J. Sci. Comput. (2015) 53: 745-768.



\bibitem{Shubin01}
M.A. Shubin  (2001) Asymptotic Behaviour of the Spectral Function. In: Pseudodifferential Operators and Spectral Theory. Springer, Berlin, Heidelberg

\bibitem{Stam03}
J. Stam, 
Flows on surfaces of arbitrary topology, 
ACM Trans. Graph. 22 (2003)  724--731.





\bibitem{Turk91}
G. Turk. 
Generating textures on arbitrary surfaces using reaction-diffusion, 
volume
25 ACM  1991.







\bibitem{Turing52}
M. A. Turing, The chemical basis of morhogenesis,Phil. Trans. R. Soc. (1952) B237,37-72. 





\end{thebibliography}
\end{document}